\documentclass[a4paper, 11pt, twoside]{article}
\usepackage{amsmath}
\usepackage{amscd}
\usepackage{amssymb}
\usepackage{amsthm}
\usepackage{graphics}
\usepackage{indentfirst}
\usepackage{latexsym}
\usepackage{amsfonts}
\usepackage{multicol}
\usepackage{stmaryrd}
\usepackage{array}
\usepackage{genyoungtab}
\usepackage{pxfonts}
\usepackage{verbatim}
\usepackage{bbm}

\makeatletter
\newenvironment{pf}[1][\proofname] {\par\pushQED{\qed}\normalfont\topsep6\p@\@plus6\p@\relax\trivlist\item[\hskip\labelsep\bfseries#1\@addpunct{.}]\ignorespaces}{\popQED\endtrivlist\@endpefalse}
\makeatother

\newtheoremstyle{mattthm}{}{}{\itshape}{}{\bfseries}{.}{ }{}

\theoremstyle{mattthm}
\newtheorem{lemma}{Lemma}[section]

\newtheorem{thm}[lemma]{Theorem}

\newtheoremstyle{mattdef}{}{}{}{}{\bfseries}{.}{ }{}

\theoremstyle{mattdef}

\begin{document}

\newenvironment{pfenum}{\begin{pf}\indent\begin{enumerate}\vspace{-\topsep}}{\qedhere\end{enumerate}\end{pf}}

\newcommand\ind{\operatorname{Ind}}
\newcommand\res{\operatorname{Res}}
\newcommand\iind{\operatorname{-Ind}}
\newcommand\ires{\operatorname{-Res}}
\newcommand\syd{\,\triangle\,}
\newcommand\mo{{-}1}
\newcommand\reg{\mathcal{R}}
\newcommand\swed[1]{\stackrel{#1}{\wed}}
\newcommand\inv[2]{\llbracket#1,#2\rrbracket}
\newcommand\one{\mathbbm{1}}
\newcommand\bsm{\begin{smallmatrix}}
\newcommand\esm{\end{smallmatrix}}
\newcommand{\rt}[1]{\rotatebox{90}{$#1$}}
\newcommand\la\lambda
\newcommand\id{\operatorname{identity}}
\newcommand{\ol}{\overline}
\newcommand{\ul}{\check}
\newcommand{\lan}{\langle}
\newcommand{\ran}{\rangle}
\newcommand\partn{\mathcal{P}}
\newcommand\calp{\mathcal{P}}
\newcommand\calf{\mathcal{F}}
\newcommand{\py}[3]{\,_{#1}{#2}_{#3}}
\newcommand{\pyy}[5]{\,_{#1}{#2}_{#3}{#4}_{#5}}
\newcommand{\thmlc}[3]{\textup{\textbf{(\!\! #1 \cite[#3]{#2})}}}
\newcommand{\thmlabel}[1]{\textup{\textbf{#1}}\ }
\newcommand{\sss}{\mathfrak{S}_}
\newcommand{\dom}{\trianglerighteqslant}
\newcommand{\doms}{\vartriangleright}
\newcommand{\ndom}{\ntrianglerighteqslant}
\newcommand{\ndoms}{\not\vartriangleright}
\newcommand{\domby}{\trianglelefteqslant}
\newcommand{\domsby}{\vartriangleleft}
\newcommand{\ndomby}{\ntrianglelefteqslant}
\newcommand{\ndomsby}{\not\vartriangleleft}
\newcommand{\subs}[1]{\subsection{#1}}
\newcommand{\subsubs}[1]{\subsubsection*{#1}}
\newcommand{\nin}{\notin}
\newcommand{\nchar}{\operatorname{char}}
\newcommand{\thmcite}[2]{\textup{\textbf{\cite[#2]{#1}}}\ }
\newcommand\zez{\mathbb{Z}/e\mathbb{Z}}
\newcommand\zepz{\mathbb{Z}/(e+1)\mathbb{Z}}
\newcommand{\bbf}{\mathbb{F}}
\newcommand{\bbg}{\mathbb{G}}
\newcommand{\bbc}{\mathbb{C}}
\newcommand{\bbn}{\mathbb{N}}
\newcommand{\bbq}{\mathbb{Q}}
\newcommand{\bbz}{\mathbb{Z}}
\newcommand\zo{\bbn_0}
\newcommand{\gs}{\geqslant}
\newcommand{\ls}{\leqslant}
\newcommand\dw{^\triangle}
\newcommand\wod{^\triangledown}
\newcommand{\hhh}{\mathcal{H}_}
\newcommand{\bbb}{\mathcal{B}_}
\newcommand{\aaa}{\mathcal{A}_}
\newcommand{\sect}[1]{\section{#1}}
\newcommand{\ff}{\mathfrak{f}}
\newcommand{\fff}{\mathfrak{F}}
\newcommand\cf{\mathcal{F}}
\newcommand\fkn{\mathfrak{n}}
\newcommand\fkp{\mathfrak{p}}
\newcommand\sx{x}
\newcommand\bra[1]{|#1\ran}
\newcommand\arb[1]{\widehat{\bra{#1}}}
\newcommand\foc[1]{\mathcal{F}_{#1}}
\newcommand{\clam}{\begin{description}\item[\hspace{\leftmargin}Claim.]}
\newcommand{\prof}{\item[\hspace{\leftmargin}Proof.]}
\newcommand{\malc}{\end{description}}
\newcommand\ppmod[1]{\ (\operatorname{mod}\ #1)}
\newcommand\wed\wedge
\newcommand\wede\barwedge
\newcommand\uu[1]{\,\begin{array}{|@{\,}c@{\,}|}\hline #1\\\hline\end{array}\,}
\newcommand{\ux}[1]{\operatorname{ht}_{#1}}
\newcommand\erim{\operatorname{rim}}
\newcommand\mire{\operatorname{rim}'}
\newcommand\mmod{\ \operatorname{Mod}}
\newcommand\cgs\succcurlyeq
\newcommand\cls\preccurlyeq
\newcommand\inc{\mathfrak{A}}
\newcommand\fsl{\mathfrak{sl}}
\newcommand\ba{\mathbf{s}}
\newcommand\ta{\tilde\ba}
\newcommand\kt[1]{|#1\rangle}
\newcommand\tk[1]{\langle#1|}
\newcommand\ket[1]{s_{#1}}
\newcommand\jn\diamond
\newcommand\UU{\mathcal{U}}
\newcommand\MM[1]{M^{\otimes#1}}
\newcommand\add{\operatorname{add}}
\newcommand\rem{\operatorname{rem}}
\newcommand\La\Lambda
\newcommand\lra\longrightarrow
\newcommand\lexg{>_{\operatorname{lex}}}
\newcommand\lexgs{\gs_{\operatorname{lex}}}
\newcommand\lexl{<_{\operatorname{lex}}}
\newcommand\tru[1]{{#1}_-}
\newcommand\ste[1]{{#1}_+}
\newcommand\out{^{\operatorname{out}}}
\newcommand\lad{\mathcal{L}}
\newcommand\hsl{\widehat{\mathfrak{sl}}}
\newcommand\fkh{\mathfrak{h}}
\newcommand\GG{H}
\newcommand\dable{restrictable}
\newcommand\infi[1]{$(\infty,#1)$-irreducible}
\newcommand\infr[1]{$(\infty,#1)$-reducible}
\newcommand\infs[1]{$#1$-signature}
\newcommand\be[2]{B^{#1}(#2)}
\newcommand\domi{dominant}
\newcommand\app{good}
\newcommand\lset[2]{\left\{\left.#1\ \right|\ #2\right\}}
\newcommand\rset[2]{\left\{#1\ \left|\ #2\right.\right\}}
%abacus drawing commands - use a smallmatrix, with \bd for a bead and \nb for a space.
\newcommand\tl{\begin{picture}(8,4)
\put(4,-1){\line(0,1){3}}
\put(4,2){\line(1,0){7}}
\end{picture}}
\newcommand\tr{\begin{picture}(8,4)
\put(4,-1){\line(0,1){3}}
\put(-3,2){\line(1,0){7}}
\end{picture}}
\newcommand\tm{\begin{picture}(8,4)
\put(4,-1){\line(0,1){3}}
\put(-3,2){\line(1,0){14}}
\end{picture}}
\newcommand{\bd}{\begin{picture}(8,6)
\put(4,-1){\line(0,1){8}}
\put(4,3){\circle*{6}}
\end{picture}}
\newcommand{\nb}{\begin{picture}(8,6)
\put(4,-1){\line(0,1){8}}
\put(3,3){\line(1,0){2}}
\end{picture}}
\newcommand{\vd}{\begin{picture}(8,10)
\put(4,5){\circle*{1}}
\put(4,2){\circle*{1}}
\put(4,8){\circle*{1}}
\end{picture}}
\newcommand{\hd}{\begin{picture}(10,6)
\put(5,3){\circle*{1}}
\put(2,3){\circle*{1}}
\put(8,3){\circle*{1}}
\end{picture}}
\newcommand{\hhd}{\begin{picture}(10,12)
\put(5,6){\circle*{1}}
\put(2,6){\circle*{1}}
\put(8,6){\circle*{1}}
\end{picture}}
\newcommand{\hbd}{\begin{picture}(10,12)
\put(-2,6){\line(1,0){14}}
\put(5,6){\circle*{10}}
\end{picture}}
\newcommand{\hnb}{\begin{picture}(10,12)
\put(-2,6){\line(1,0){14}}
\put(5,4.5){\line(0,1){3}}
\end{picture}}

%Topmatter
\title{The exact spread of M$_{23}$ is 8064}
\author{Ben Fairbairn\\\normalsize Departmento de Matem\'{a}ticas, Universidad de los Andes,\\ Carrera 1 No 18A-12, Bogot\'{a}, Colombia\\\texttt{\normalsize bt.fairbairn20@uniandes.edu.co}}
\date{}
\maketitle
\markboth{Ben Fairbairn}{The exact spread of M$_{23}$ is 8064}
\pagestyle{myheadings}

\begin{abstract}
We show that if $\{x_1,\ldots,x_{8064}\}\subset$ M$_{23}^{\#}$ then there exists an element $y\in$ M$_{23}$ such that $\langle y,x_i\rangle=$ M$_{23}$ for every $1\leq i\leq8064$ and that no larger set of distinct elements of M$_{23}$ has this property. 
\end{abstract}

\section{Introduction}

Let $G$ be a group. We say that $G$ has \emph{spread r} if for any
set of distinct non-trivial elements $X:=\{x_1,\ldots, x_r\}\subset
G^{\#}$ there exists an element $y\in G$ with the property that
$\langle x_i,y\rangle=G$ for every $1\leq i\leq r$. We say that $y$ is a \emph{mate} to $X$. We say $G$ has
\emph{exact spread} $r:=s(G)$ if $G$ has spread $r$ but not $r+1$.

This concept was first introduced by Brenner and Wiegold in \cite{BrennerWiegold}, extending earlier work of Binder in \cite{Binder}. The concept of spread may be thought of as a generalization of the idea of 3/2-generation. 

Since this note is not intended for publication, we omit background and motivation, but the interested reader may wish to consult \cite{BradleyHolmes,BradleyMoori,BreuerGuralnicKantor,F,GaniefMoori,GuralnicShalev,Woldar}. The reader should at least note that there has been considerable interest in determining, or at least bounding, the exact spreads of the finite simple groups and in particular the exact spreads of the sporadic simple groups. The precise value of the exact spread of a finite simple group is known in very few cases. Here we prove the following.

\begin{thm}
$s(\mbox{M}_{23})=8064$
\end{thm}

\section{Proof of Theorem 1}

\begin{proof}
First note that from \cite[Table 1]{F} (itself a heavily corrected version of \cite[Table 1]{BradleyHolmes}) we have that $s(\mbox{M}_{23})\leq8064$ and so it is sufficient to show that any set of 8064 elements from M$_{23}^{\#}$ has a mate. Let $X\subset$ M$_{23}$ be a set of 8064 distinct elements. 

Now, from the maximal subgroups of M$_{23}$ listed in \cite[p.71]{ATLAS} we see that an element of order 23 is contained in only one maximal subgroup - a copy of the Frobenius group 23:11. Since the normalizer in M$_{23}$ of a cyclic subgroup of order 11 is a Frobenius group 11:5, each element of order 11 is contained in five disctinct copies of 23:11. It follows that the only way $X$ can avoid having a mate of order 23 is if $X$ consists of a well chosen configutation of elements of order 11.

Finally, note that of the maximal subgroups of M$_{23}$ containing elements of order 11, which are each isomorphic to one of M$_{22}$, M$_{11}$ or 23:11, none contain elements of order 14. It follows that even if $X$ only contains elements of order 11, $X$ must have a mate and so $s(\mbox{M}_{23})\geq8064$.
\end{proof}

\end{document}